\documentclass{article}

\usepackage[repositionpullbacks,midshaft,PostScript=dvips,nohug,scriptlabels,size=2em]{diagrams}

\usepackage{latexsym}
\usepackage{amssymb}
\usepackage{mathrsfs}
\usepackage{graphics}

\newcommand{\bref}[1]{(\ref{#1})}

\newcommand{\cat}[1]{\mathscr{#1}}
\newcommand{\fcat}[1]{\mathbf{#1}}

\newcommand{\such}{\mid}
\newcommand{\iso}{\cong}
\newcommand{\of}{\,\raisebox{0.08ex}{\ensuremath{\scriptstyle\circ}}\,}
\newcommand{\ladj}{\dashv}
\newcommand{\id}{\mathrm{id}}
\newcommand{\Aut}{\mathrm{Aut}}
\newcommand{\End}{\mathrm{End}}

\newcommand{\Set}{\fcat{Set}}
\newcommand{\Cat}{\fcat{Cat}}
\newcommand{\Multicat}{\fcat{Multicat}}
\newcommand{\Operad}{\fcat{Operad}}
\newcommand{\Span}{\fcat{Span}}
\newcommand{\MonCat}{\fcat{MonCat}}
\newcommand{\MonGpd}{\fcat{MonGpd}}
\newcommand{\Gpd}{\fcat{Gpd}}

\newcommand{\nat}{\mathbb{N}}	

\newcommand{\pull}{\mathrm{pb}}
\newcommand{\card}[1]{| #1 |}
\newcommand{\bintr}[1]{\fcat{Tr}_{#1}}

\newcommand{\slsh}{/\linebreak[0]}

\newcommand{\dt}{.\linebreak[0]}
\newcommand{\cln}{:\linebreak[0]}

\newtheorem{thm}{Theorem}[section]
\newtheorem{propn}[thm]{Proposition}
\newtheorem{lemma}[thm]{Lemma}

\newcommand{\demph}[1]{\textbf{\textup{#1}}} 
\newcommand{\latin}[1]{#1}
\newcommand{\smallheading}[1]{\subsubsection*{#1}}
\newenvironment{univ}{\begin{quote}}{\end{quote}}
\newcommand{\littlespace}{\,}

\newcommand{\done}{\hfill\ensuremath{\Box}}
\newenvironment{prooflike}[1]{\begin{trivlist}\item\textbf{#1}\ }
{\end{trivlist}}
\newenvironment{proof}{\begin{prooflike}{Proof}}{\end{prooflike}}

\newcommand{\go}{\rTo\linebreak[0]}
\newcommand{\goby}[1]{\rTo^{#1}\linebreak[0]}
\newcommand{\goesto}{\,\longmapsto\,}
\newcommand{\oppair}[2]{\pile{\rTo^{\scriptstyle #1}\\ 
\lTo_{\scriptstyle #2}}} 

\newcommand{\adjnpair}{\pile{\rTo\\ \lTo}}
\newcommand{\adjntriple}{\pile{\rTo\\ \lTo\\ \rTo}}
\newarrow{Goesto}|---{->}
\newarrow{Incl}C--->
\newarrowmiddle{iso}\sim\sim\sim\sim

\newcommand{\ternarytree}{%
\setlength{\unitlength}{1ex}
\begin{picture}(2.1,2.7)(0,0)
\cell{1.3}{0}{b}{\bf\sf Y}
\cell{0.63}{1.5}{b}{\bf\sf v}
\end{picture}}
\newcommand{\binarytree}{\mbox{\bf\sf Y}}
\newcommand{\utree}{\mbox{\bf\sf l}}

\newcommand{\cell}[4]{\put(#1,#2){\makebox(0,0)[#3]{\ensuremath{#4}}}}
\newcommand{\zmark}{\scriptstyle{\bullet}}

\newcommand{\tusual}[1]{%
\begin{picture}(4,4)(-2,-2)
\cell{-0.2}{0}{c}{#1}%
\put(-2,-2){\line(0,1){4}}%
\put(-2,2){\line(2,-1){4}}%
\put(2,0){\line(-2,-1){4}}%
\end{picture}}

\newcommand{\toutputrgt}[1]{%
\begin{picture}(1,0)(-1,0)
\cell{0.4}{0}{l}{#1}%
\put(0,0){\line(-1,0){1}}%
\end{picture}}

\newcommand{\tinputlft}[1]{%
\begin{picture}(1,0)(0,0)
\cell{-0.4}{0}{r}{#1}%
\put(0,0){\line(1,0){1}}
\end{picture}}

\newcommand{\tinputslft}[2]{%
\begin{picture}(1.4,3.8)(-0.4,-1.9)
\cell{0}{1.5}{l}{\tinputlft{#1}}%
\cell{0}{-1.5}{l}{\tinputlft{#2}}%
\cell{0.2}{0.3}{c}{\vdots}%
\end{picture}}

\title{An abstract characterization of\\
Thompson's group $F$}
\author{%
Marcelo Fiore%
\thanks{Computer Laboratory, University of Cambridge, UK;
Marcelo.Fiore@cl.cam.ac.uk.  
Partially supported by an EPSRC Advanced Research Fellowship.} 
\qquad
Tom Leinster%
\thanks{Department of Mathematics, University of Glasgow, UK; 
T.Leinster@maths.gla.ac.uk.
Partially supported by a Nuffield Foundation award NUF-NAL 04 and an EPSRC
Advanced Research Fellowship.}
}
\date{}

\begin{document}

\sloppy

\maketitle

\begin{abstract}
We show that Thompson's group $F$ is the symmetry group of the `generic
idempotent'.  That is, take the monoidal category freely generated by an
object $A$ and an isomorphism $A \otimes A \go A$; then $F$ is the group of
automorphisms of $A$.
\end{abstract}

\section{Introduction}
\label{sec:intro}

Our purpose in this paper is to clarify an idea concerning Richard Thompson's
group $F$: that it is, in a suitable sense, the automorphism group of some
object known only to be isomorphic to a combination of two copies of itself.
This general idea has been known for some years, but it does not seem to have
been observed until now that it can be formalized very succinctly.  We prove
that $F$ can be defined as follows.  Take the monoidal category freely
generated by an object $A$ and an isomorphism $A \otimes A \go A$; then $F$ is
the group of automorphisms of $A$.  This result first appeared in
our 2005 preprint~\cite{FL}.

Our characterization is distinct from some superficially similar older
characterizations.  In particular, it is distinct from Higman's
characterization of Thompson's group $V$ as the automorphism group of a
certain free algebra, and of $F$ as the subgroup consisting of the
`order-preserving' automorphisms~\cite{Hig,Bro,CFP}.  It is also distinct from
Freyd and Heller's characterization of $F$ via conjugacy
idempotents~\cite{FH}.  We do not know of any direct way to deduce our
characterization from these older ones, or \latin{vice versa}.  

Intuitively, our result means the following.  Suppose that we are handed a
mathematical object and told only that it is isomorphic to two copies of
itself glued together.  We do not know what kind of object it is, nor do we
know what `gluing' means except that it is some kind of associative operation.
On the basis of this information, what automorphisms does our object have?
Our result gives the answer: the elements of $F$.

Our description of $F$ is not only conceptually simple, but is also a member
of a well-established family: many entities of interest can be described via
free categories with structure.  For example, the braided monoidal category
freely generated by one object is the sequence $(B_n)_{n\geq 0}$ of Artin
braid groups~\cite{JS,Mac}.  The monoidal category freely generated by a
monoid consists of the finite ordinals; in other words, it is the augmented
simplex category~\cite{Lawv,Mac}.  The symmetric monoidal category freely
generated by a commutative monoid consists of the finite cardinals.  The
symmetric monoidal category freely generated by a commutative Frobenius
algebra consists of $1$-dimensional smooth oriented manifolds and
diffeomorphism classes of $2$-dimensional cobordisms.  (This last example is a
strong form of the equivalence between commutative Frobenius algebras and
$2$-dimensional topological quantum field theories~\cite{Dij};
see~\cite{Kock}, for instance.)  In this vein, our result can be expressed as
follows: the monoidal category freely generated by an object $A$ and an
isomorphism $A \otimes A \go A$ is equivalent to the groupoid $1 \amalg F$,
where $1$ is the trivial group and $\amalg$ is coproduct of groupoids.

Our result is this:
\begin{thm}	\label{thm:main}
Let $\cat{A}$ be the monoidal category freely generated by an idempotent
object $(A, \alpha)$.  Then $\Aut_{\cat{A}}(A)$ is isomorphic to Thompson's
group $F$.
\end{thm}
To make this paper accessible to as wide a readership as possible, we give the
definition of Thompson's group and explain the categorical language used in
the statement of this theorem~(\S\ref{sec:term}).  (The only new piece of
terminology is `idempotent object', which means an object $A$ together with an
isomorphism $\alpha: A \otimes A \go A$.)  But first we discuss earlier
characterizations of Thompson's group.

\smallheading{Related work}

Almost as soon as Thompson introduced the group now called $F$, it began to be
understood that $F$ was in some sense the automorphism group of an object
known only to be isomorphic to a combination of two copies of itself.  This
intuition is so crucial that it has been formalized in several ways, of which
ours is one.

An early such formalization, due to Thompson and Higman, was as follows.
A \demph{J\'onsson--Tarski algebra}~\cite{JT}, or \demph{Cantor
algebra}, is a set $A$ equipped with a bijection $A \times A \go A$.
Thompson's group $V$ is the automorphism group of the free J\'onsson--Tarski
algebra on one generator~\cite{Hig,CFP}.  Thompson's group $F$ is the subgroup
consisting of those automorphisms that are `order-preserving' in a suitable
sense \cite{Bro,CFP}.
There is a clear resemblance between these descriptions and ours.  However, we
know of no direct or simple way to deduce our description of $F$ from the
earlier one (or indeed the converse).  

There is a sense in which our description of $F$ is more direct.  Whenever one
works with sets and their cartesian products, one automatically introduces a
symmetry in the form of the natural isomorphism $X \times Y \go Y \times X$
for sets $X$ and $Y$.  In particular, for sets $X$, there is a nontrivial
natural automorphism of $X \times X$.  In Thompson and Higman's description,
symmetry is first created (by working with sets) and then destroyed (by
restricting to the order-preserving automorphisms).  In our approach symmetry
is avoided entirely, by working from the start not with sets, but with objects
of a (non-symmetric) monoidal category.

This is also what makes it possible to characterize $F$ as the \emph{full}
automorphism group of some algebraic structure, rather than just a subgroup.
As far as we know, this is the first such characterization.

Among all the results related to ours, the closest is probably a theorem of
Guba and Sapir~\cite{GS}.  Given any presentation of a monoid, they define
what they call its Squier complex, a 2-dimensional complex whose
connected-components are the elements of the monoid.  Every element of the
monoid therefore gives rise to a `diagram group', the fundamental group of the
corresponding component.  They show that the diagram group of the presentation
$\langle x \mid x^2 = x \rangle$ at the element $x$ is $F$.  The connection
between their result and ours can be summarized as follows.  First, the Squier
complex of this presentation is (up to homotopy) the $2$-skeleton of the
classifying space of the monoidal category freely generated by an idempotent
object $(A, \alpha)$.  (For explanation of the latter phrase,
see~\S\ref{sec:term}; for classifying spaces of categories, see~\cite{Seg},
for instance.)  Then, the generator $x$ determines a point of the Squier
complex, the object $A$ determines a point of the classifying space, and these
two points correspond to one another under the homotopy equivalence.  Hence
the fundamental group at $x$ is the automorphism group of $A$.  In this way,
their result can be deduced from ours and \latin{vice versa}.

Some more distant relatives are the results of Brin~\cite{Brin2}, Dehornoy
\cite{De1,De2}, and, ultimately, McKenzie and Thompson~\cite{MT}.  In the
context of semigroup theory, our work has connections with recent work of
Lawson~\cite{Laws,LawsCST}. 

All of these results express how $F$ arises naturally from two very primitive
notions: binary operation and associativity.  An advantage of our approach is
that it makes this idea precise using only standard categorical language,
where other approaches have used language invented more or less specifically
for the occasion.

A further advantage is that Thompson's group $V$, and even
higher-dimensional versions of it, have similar characterizations: for $V$,
just replace `monoidal category' by `symmetric monoidal category', or
equally `finite-product category'.  We do not know whether there is such a
characterization of Thompson's group $T$; using braided monoidal categories
gives not $T$, but the braided version of $V$ defined in~\cite{Brin1}.
Also, given any $n \geq 2$, if we take the monoidal category freely
generated by an object $A$ and an isomorphism $A^{\otimes n} \go A$ then
the automorphism group of $A^{\otimes r}$ is canonically isomorphic to the
generalized Thompson group $F_{n, r}$ of Brown~\cite{Bro}.  

Freyd and Heller also gave a short categorical definition of $F$, different
from ours: it is the initial object in the category of groups equipped with
a conjugacy-idempotent endomorphism~\cite{FH}.  Again, there is a striking
resemblance between this description and ours; but again, no one (to our
knowledge) has been able to find a direct deduction of one from the other.  

The category of forests and the free groupoid on it, which appear
in~\S\ref{sec:proof} below, have been considered independently by
Belk~\cite{Belk}.

We work throughout with \emph{strict} monoidal categories.  (See below for
definitions.)  However, the non-strict monoidal category freely generated by
an idempotent object $(A', \alpha')$ is monoidally equivalent to the strict
one, and in particular, the automorphism group of $A'$ is $F$.  So, for
instance, there is an induced homomorphism from $F$ to the automorphism group
of the free J\'onsson--Tarski algebra on one generator.  We conjecture that
this homomorphism is injective and that its image consists of the
order-preserving automorphisms.

\smallheading{}

In~\S\ref{sec:term} we explain all of the terminology used in the statement of
Theorem~\ref{thm:main}.  The theorem is proved in~\S\ref{sec:proof}.  Our
proof involves almost no calculation, but does use some further concepts from
category theory, reviewed in the Appendix.  (`\emph{Il faut triompher par la
pens\'ee et non par le calcul'}%
\footnote{One must prevail by thought, not by calculation.}%
---Poincar\'e.)  

Some readers may feel that the language used in the statement of the theorem
represents quite enough category theory for their taste, even without 
the further categorical concepts used in the proof.  For them we sketch, at
the end of~\S\ref{sec:term}, an alternative proof, favouring explicit
calculation over conceptual argument.

The novelty of this work lies almost entirely in~\S\ref{sec:proof} and in the
way in which the categorical and algebraic structures are brought together.
In particular, the categorical language explained in~\S\ref{sec:term} is
absolutely standard; and while not everything in the Appendix is quite so
well known, none of it is by any means new.

\section{Terminology}
\label{sec:term}
Here we explain the terminology in the statement of
Theorem~\ref{thm:main}.  Further information on Thompson's group can be found
in~\cite{CFP}; for more on the categorical language, see~\cite{Mac}.  We then
sketch a calculational proof of Theorem~\ref{thm:main}, requiring no further
categorical concepts.

\begin{figure}
\centering
\setlength{\unitlength}{4ex}
\begin{picture}(8,3)(-4,-1.5)
\thinlines
\put(-6,-0.5){\line(6,-1){6}}
\put(-6,0.5){\line(6,1){6}}
\put(-3,0){\line(6,-1){3}}
\put(6,-0.5){\line(-6,-1){6}}
\put(6,0.5){\line(-6,1){6}}
\put(3,0){\line(-6,1){3}}
\thicklines
\put(-6,-0.5){\line(0,1){1}}
\put(-3,-1){\line(0,1){2}}
\put(0,-1.5){\line(0,1){3}}
\put(3,-1){\line(0,1){2}}
\put(6,-0.5){\line(0,1){1}}
\cell{-6.1}{0.5}{r}{0}
\cell{-6.1}{-0.5}{r}{1}
\cell{-3.1}{1}{r}{0}
\cell{-3.1}{0}{r}{1}
\cell{-3.1}{-1}{r}{2}
\cell{-0.1}{1.5}{r}{0}
\cell{-0.1}{0.5}{r}{1}
\cell{-0.1}{-0.5}{r}{2}
\cell{-0.1}{-1.5}{r}{3}
\cell{3.1}{1}{l}{0}
\cell{3.1}{0}{l}{1}
\cell{3.1}{-1}{l}{2}
\cell{6.1}{0.5}{l}{0}
\cell{6.1}{-0.5}{l}{1}
\cell{-4.5}{0}{c}{\alpha^{-1}}
\cell{-1.5}{0.5}{c}{\alpha^{-1}}
\cell{-1.5}{-0.75}{c}{\id}
\cell{1.5}{0.75}{c}{\id}
\cell{1.5}{-0.5}{c}{\alpha}
\cell{4.5}{0}{c}{\alpha}
\end{picture}
\caption{Decomposition of an element of $F$}
\label{fig:decomp}
\end{figure}
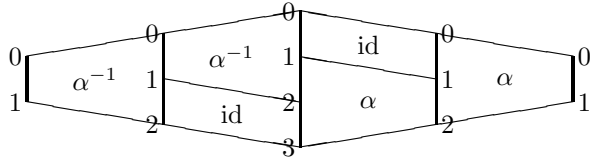

\smallheading{Thompson's group $F$}  
In the 1960s Richard Thompson discovered three groups, now called $F$, $T$
and $V$, with remarkable properties.  The group $F$, in particular, is
one of those mathematical objects that appears in many diverse contexts
and has been rediscovered repeatedly.  One definition of $F$ is
that it consists of all bijections $f: [0, 1] \go [0, 1]$ satisfying
\begin{enumerate}
\item \label{item:F-defn-one}
$f$ is piecewise linear (with only finitely many pieces)
\item
the slope (gradient) of each piece is an integer power of $2$
\item \label{item:F-defn-three}
the coordinates of the endpoints of each piece are dyadic rationals.
\end{enumerate}
For example, the $3$-piece linear function $f$ satisfying $f(0) =
0$, $f(1/4) = 1/2$, $f(1/2) = 3/4$ and $f(1) = 1$ is an element
of $F$.  In a sense that will be made precise, every element of
$F$ can be built from copies of the halving isomorphism $\alpha:
[0, 2] \go [0, 1]$ and its inverse; this is shown for our example
$f$ in Figure~\ref{fig:decomp}.
So if all we knew about $[0, 1]$ was that it was isomorphic 
to two copies of itself glued together, $F$ would be
the group of \emph{all} automorphisms of $[0, 1]$.  This is the
spirit of our result.

\begin{figure}
\centering
\setlength{\unitlength}{3ex}
\begin{picture}(11,4.5)(0,-2.5)
\thicklines
\put(0,0){\line(1,0){1}}
\put(1,0){\line(2,1){4}}
\put(1,0){\line(2,-1){4}}
\put(3,1){\line(2,-1){2}}
\cell{1}{0}{c}{\zmark}
\cell{3}{1}{c}{\zmark}
\put(11,0){\line(-1,0){1}}
\put(10,0){\line(-2,1){4}}
\put(10,0){\line(-2,-1){4}}
\put(8,-1){\line(-2,1){2}}
\cell{10}{0}{c}{\zmark}
\cell{8}{-1}{c}{\zmark}
\cell{3}{-2.5}{b}{\tau'}
\cell{8}{-2.5}{b}{\tau}
\end{picture}
\caption{A representative $(\tau, \tau')$ of an element of $F$}
\label{fig:twotrees}
\end{figure}
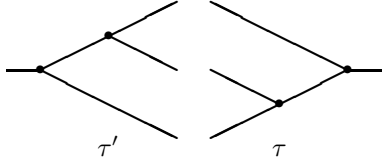

For the proof we will need an alternative, more combinatorial definition of
$F$.  In what follows, \demph{tree} will mean finite, rooted, planar
tree, and a tree is \demph{binary} if precisely two branches grow out of
each vertex.  Figure~\ref{fig:twotrees} shows a pair of binary trees.  Except
where mentioned, `tree' will mean `binary tree'.

For $n \in \nat = \{0, 1, 2, \ldots\}$, write $\bintr{n}$ for the set of
$n$-leafed trees.  There are no $0$-leafed trees, and there is just
one $1$-leafed tree: the trivial tree $\utree$ with no vertices at all.  A
non-trivial tree consists of two smaller trees joined at the root, so the sets
$\bintr{n}$ can be defined inductively by
\[
\bintr{0} = \emptyset,
\quad
\bintr{1} = \{ \utree \},
\quad
\bintr{n} = \coprod_{k + m = n} \bintr{k} \times \bintr{m}
\quad
(n \geq 2).  
\]
By a \demph{subtree} of a tree we mean a subtree sharing the same root.  For
example, the tree $\ternarytree$ has exactly three subtrees: itself, the unique
two-leafed tree $\binarytree$, and the one-leafed tree $\utree$. 

Given $n \geq i \geq 1$, we can join to the $i$th leaf of any $n$-leafed tree
$\tau$ a copy of the two-leafed tree $\binarytree$, thus forming an $(n +
1)$-leafed tree $\omega^n_i(\tau)$.  This defines a map $\omega^n_i: \bintr{n}
\go \bintr{n+1}$.  Whenever $\tau$ is a subtree of a tree $\rho$, there is a
finite sequence $\omega^{n_1}_{i_1}, \ldots, \omega^{n_r}_{i_r}$ of maps such
that
\[
\rho
=
\omega^{n_r}_{i_r} \cdots \omega^{n_1}_{i_1} (\tau).
\]
Moreover, for any two trees $\sigma$ and $\tau$, there is a smallest
tree containing both as subtrees.  This
can be obtained by superimposing the pictures of $\sigma$ and $\tau$. 

The following alternative definition of $F$ is given in~\cite[\S 2]{CFP} and
in~\cite[1.2]{Belk}.  Elements of $F$ are equivalence classes of pairs $(\tau,
\tau')$ of trees with the same number of leaves, where the equivalence
relation is generated by identifying $(\tau, \tau')$ with $(\omega^n_i (\tau),
\omega^n_i (\tau'))$ whenever $\tau, \tau' \in \bintr{n}$ and $1 \leq i \leq
n$.  Write $[\tau, \tau']$ for the equivalence class of a pair $(\tau,
\tau')$.

Under this definition, the element of $F$ shown in Figure~\ref{fig:decomp} is
the same as the element $[\tau, \tau']$ shown in~Figure~\ref{fig:twotrees}.
In general, $[\tau, \tau']$ can be read as `expand according to $\tau'$
then contract according to $\tau$'.  (The order is reversed to agree with
the convention of writing maps on the left.)

With this in mind, it is clear what the product (composite) $[\tau,
\tau']\, [\sigma, \sigma']$ must be when $\tau' = \sigma$: simply $[\tau,
\sigma']$.  In general, there is a tree containing both
$\tau'$ and $\sigma$ as subtrees, so there are maps $\omega^{n_q}_{i_q},
\omega^{m_q}_{j_q}$ for which
\[
\omega^{n_r}_{i_r} \cdots \omega^{n_1}_{i_1} (\tau')
=
\omega^{m_s}_{j_s} \cdots \omega^{m_1}_{j_1} (\sigma),
\]
and then---inevitably---
\[
[\tau, \tau']\, [\sigma, \sigma']
=
[ \omega^{n_r}_{i_r} \cdots \omega^{n_1}_{i_1} (\tau),\littlespace
\omega^{m_s}_{j_s} \cdots \omega^{m_1}_{j_1} (\sigma') ].
\]

\smallheading{Monoidal categories} 
A monoid is a set $S$ equipped with a function $S \times S \go S$ and an
element $1 \in S$ obeying associativity and unit laws.  Similarly, a
\demph{monoidal category} is a category $\cat{M}$ equipped with a functor
$\cat{M} \times \cat{M} \go \cat{M}$ and an object $I \in \cat{M}$ obeying
associativity and unit laws.  Explicitly, this means that to each pair $(M,
N)$ of objects of $\cat{M}$ there is assigned an object $M \otimes N$, and
to each pair
\[
\left(
M \goby{\phi} M',
\ 
N \goby{\psi} N'
\right)
\]
of maps in $\cat{M}$ there is assigned a map $M \otimes N \goby{\phi
\otimes \psi} M' \otimes N'$.  Functoriality amounts to the equations
\[
(\phi' \of \phi) \otimes (\psi' \of \psi) 
=
(\phi' \otimes \psi') \of (\phi \otimes \psi),
\quad
1_M \otimes 1_N 
=
1_{M \otimes N},
\]
and the associativity and unit laws apply to maps as well as objects: $(\phi
\otimes \psi) \otimes \chi = \phi \otimes (\psi \otimes \chi)$, etc.  A
\demph{monoidal functor} is a functor $G$ between monoidal categories that
preserves the tensor and unit: $G(M \otimes N) = G(M) \otimes G(N)$, etc.

For example, a monoidal category in which the only maps are identities is
simply a monoid.  The monoidal category $\fcat{FinOrd}$ of finite ordinals has
as objects the natural numbers; a map $m \go n$ is an order-preserving
function $\{ 0, \ldots, m-1 \} \go \{ 0, \ldots, n-1\}$; the tensor product is
given on objects by addition and on maps by juxtaposition; the unit object is
$0$.

The monoidal categories and functors considered in this paper are properly
called \emph{strict} monoidal.  The more general notion of monoidal category
includes such examples as the category of abelian groups, in which the tensor
product is only associative and unital up to (suitably coherent) isomorphism.

\smallheading{Freely generated} 
We defined $\cat{A}$ as the `monoidal category freely generated by an
idempotent object $(A, \alpha)$'.  Such use of language is standard in
category theory, and extends the familiar notion of free structure in
algebra.  We now explain what it means.

Informally, it means that $\cat{A}$ is constructed by starting with an object
$A$ and an isomorphism $\alpha: A \otimes A \go A$, then adjoining whatever
other objects and maps must be present in order for $\cat{A}$ to be a monoidal
category.  The only equations that hold are those that are forced to hold by
the axioms for a monoidal category.  Thus, $\cat{A}$ has an object $A$, so it
also has an object $A^{\otimes n} = A \otimes \cdots \otimes A$ for each $n
\geq 0$ (with $A^{\otimes 0} = I$).  The maps are built up from $\alpha$ by
taking composites, identities, inverses and tensor products: for instance,
there is a map $A \go A$ given as the composite
\[
A				\goby{\alpha^{-1}} 
A \otimes A			\goby{\alpha^{-1} \otimes 1_A}
A \otimes A \otimes A		\goby{1_A \otimes \alpha}
A \otimes A			\goby{\alpha}
A.
\]
(Compare Figures~\ref{fig:decomp} and~\ref{fig:twotrees}.)

Precisely, an \demph{idempotent object} in a monoidal category $\cat{M}$ is an
object $M \in \cat{M}$ together with an isomorphism $\mu: M \otimes M \go M$.
(For example, an idempotent object in the monoidal category of sets, where
$\otimes$ is cartesian product, is a J\'onsson--Tarski algebra.)  A
\demph{monoidal category freely generated by an idempotent object} is a
monoidal category $\cat{A}$ together with an idempotent object $(A, \alpha)$
in $\cat{A}$, satisfying the following universal property:
\begin{univ}
\label{p:univ-prop}
for any monoidal category $\cat{M}$ and idempotent object $(M, \mu)$ in
$\cat{M}$, there is a unique monoidal functor $G: \cat{A} \go \cat{M}$ such
that $G(A) = M$ and $G(\alpha) = \mu$.  
\end{univ}
The universal property determines $(\cat{A}, A, \alpha)$ uniquely, up to
isomorphism.  That such an $(\cat{A}, A, \alpha)$ exists at all is true for
quite general categorical reasons, although in fact we will construct it
explicitly.  We call $(A, \alpha)$ the \demph{generic idempotent object}.

Specifying a monoidal category in this fashion is closely analogous to what
one does in algebra when specifying a group, monoid, etc.\ by a
presentation.  Suppose, say, that we define a monoid $E$ by the presentation
$E = \langle e \such e^2 = e \rangle$.  Informally, this means that $E$ is
constructed by starting with an element $e$, then adjoining whatever other
elements must be present in order for $E$ to be a monoid, then imposing only
those equations that are forced to hold by $e^2 = e$ and the axioms for a
monoid.  (Of course, for this particular presentation it is very easy to
describe $E$ explicitly, but for other presentations it is not.)  Precisely,
it means that $E$ is a monoid equipped with an idempotent element $e$ and
satisfying the following universal property:
\begin{univ}
for any monoid $X$ and idempotent element $x \in X$, there is a unique
monoid homomorphism $g: E \go X$ such that $g(e) = x$.
\end{univ}
We might call $E$ the `monoid freely generated by an idempotent element', and
$e$ the `generic idempotent element', since it is idempotent and satisfies no
further equations.

Our definition of $\cat{A}$ can be regarded as a categorification of the
definition of $E$.  Monoids have become monoidal categories, elements have
become objects, monoid homomorphisms have become monoidal functors, and
equations (such as $e^2 = e$) have become isomorphisms (such as $\alpha: A
\otimes A \go A$).  

For any monoid $X$, there is a natural one-to-one
correspondence between idempotent elements of $X$ and homomorphisms $E \go
X$.  Similarly, for any monoidal category $\cat{M}$, there is a natural
one-to-one correspondence between idempotent objects in $\cat{M}$ and monoidal
functors $\cat{A} \go \cat{M}$.

\smallheading{Automorphism group} 
Any object $X$ of any category $\cat{X}$ has an \demph{automorphism group}
$\Aut_{\cat{X}}(X)$.  Its elements are the automorphisms of $X$, that is, the
isomorphisms $X \go X$ in $\cat{X}$.  The group structure is given by
composition.

\smallheading{}

This completes the explanation of the language used in Theorem~\ref{thm:main}.
We are now in a position to sketch a proof of the theorem based on explicit
calculation, which we do for the reasons stated at the end of the
Introduction.

Let $\cat{A}$ be the category whose objects are the natural numbers and whose
maps $m \go n$ are the bijections $f: [0, m] \go [0, n]$ satisfying
conditions~(\ref{item:F-defn-one})--(\ref{item:F-defn-three}) in the
definition of $F$.  Then $\cat{A}$ has a monoidal structure given on objects
by addition and on maps by juxtaposition, and there is an isomorphism $\alpha:
1 \otimes 1 = 2 \go 1$ given by division by $2$.  We have $F =
\Aut_{\cat{A}}(1)$ by definition, so our task is to show that $(\cat{A}, 1,
\alpha)$ has the universal property stated above.

To do this, first consider trees (binary, as usual).  Take a monoidal category
$\cat{M}$ and an idempotent object $(M, \mu)$ in $\cat{M}$.  Then any
$n$-leafed tree $\tau$ gives rise to an isomorphism $\mu_\tau: M^{\otimes n}
\go M$; for instance, if $\tau = \binarytree$ then $\mu_\tau = \mu$, and if
$\tau = \ternarytree$ then $\mu_\tau$ is the composite
\[
M \otimes M \otimes M   \goby{\mu \otimes 1}
M \otimes M                     \goby{\mu}
M.
\]
More generally, define a \demph{forest} to be a finite sequence
$(\tau_1, \ldots, \tau_k)$ of trees ($k \geq 0$), and let us say that
this forest has $n$ \demph{leaves} and $k$ \demph{roots}, where $n$ is the sum
of the numbers of leaves of $\tau_1, \ldots, \tau_k$.  Any forest $T =
(\tau_1, \ldots, \tau_k)$ with $n$ leaves and $k$ roots induces an isomorphism
\[
\mu_T
=
\mu_{\tau_1} \otimes\cdots\otimes \mu_{\tau_k}:
M^{\otimes n} \go M^{\otimes k}.
\]

Now, it can be shown that any map $\phi: m \go n$ in $\cat{A}$ 
factorizes as
\begin{equation}	\label{eq:factorization}
\phi 
=
\left(
m \goby{\alpha_S^{-1}}
p \goby{\alpha_T}
n
\right)
\end{equation}
for some $p \in \nat$ and forests $S$ and $T$.  (The method is
given in~\cite[\S 2]{CFP}.)  It can also be shown that any monoidal
functor $G: \cat{A} \go \cat{M}$ satisfying $G(1) = M$ and $G(\alpha) =
\mu$ must also satisfy
\begin{equation}        \label{eq:image-map}
G(\phi)
=
\left(
M^{\otimes m} \goby{\mu_S^{-1}}
M^{\otimes p} \goby{\mu_T}
M^{\otimes n}
\right).
\end{equation}
Although $\phi$ may have many factorizations of the
form~(\ref{eq:factorization}), further calculations show that the right-hand
side of~(\ref{eq:image-map}) is independent of the factorization chosen.
Further calculations still show that the $G$ thus defined is a functor, and
monoidal.  The result follows.

\section{Proof of the Theorem}
\label{sec:proof}

In this section we give a conceptual proof of Theorem~\ref{thm:main}.

\begin{figure}
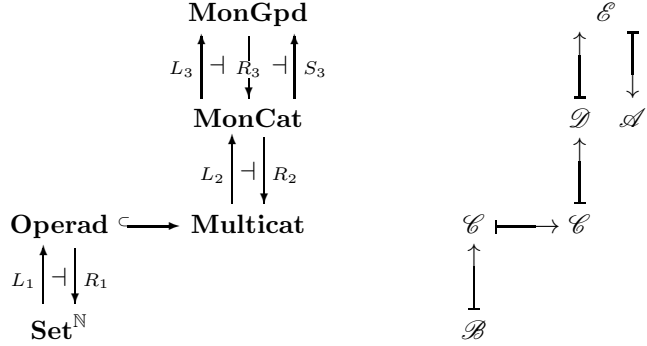

\[
\begin{diagram}
		&	&\MonGpd	\\
		&	&\uTo<{L_3} \ladj \dTo~{R_3} \ladj \uTo>{S_3}	\\
		&	&\MonCat	\\
		&	&\uTo<{L_2} \ladj \dTo>{R_2}	\\
\Operad		&\rIncl	&\Multicat	\\
\uTo<{L_1} \ladj \dTo>{R_1}&	&	\\
\Set^\nat	&	&		\\
\end{diagram}
\hspace*{5em}
\begin{diagram}
		&	&\              &\makebox[0em]{\hspace*{-2em}\ensuremath{\cat{E}}}
					&\       	\\
		&	&\uGoesto	&\dGoesto	\\
		&	&\cat{D}	&\cat{A}	\\
		&	&\uGoesto	&		\\
\cat{C}		&\rGoesto&\cat{C}	&		\\
\uGoesto	&	&		&		\\
\cat{B}		&	&		&		\\
\end{diagram}
\]
\caption{Steps in the proof}
\label{fig:steps}
\end{figure}

To do this, we construct the monoidal category $\cat{A}$ freely generated by
an idempotent object $(A, \alpha)$.  The strategy is to start with a very
simple object $\cat{B}$ and apply several left adjoints in succession
(Figure~\ref{fig:steps}).
On the one hand, this abstract construction makes the universal property of
$\cat{A}$ automatic.  On the other, each step of the construction can be
described explicitly, so it will be transparent that $\Aut_{\cat{A}}(A) \iso
F$.

On the left of Figure~\ref{fig:steps}, we have the category $\Set^\nat$ of
`signatures' and the categories of operads, multicategories, monoidal
categories and monoidal groupoids, all non-symmetric.  The functors $R_i$ are
the evident forgetful functors; they have adjoints $L_i$ and $S_3$ as shown.
Definitions, and descriptions of these adjoint functors, are given in
the Appendix.  On the right of Figure~\ref{fig:steps}, the signature $\cat{B}$
consists of a single binary operation: $\card{\cat{B}_2} = 1$ and
$\card{\cat{B}_n} = 0$ for $n \neq 2$.  Then $\cat{C} = L_1(\cat{B})$, etc.;
thus, the monoidal category $\cat{A}$ is defined by
\[
\cat{A} 
= 
R_3 L_3 L_2 L_1 (\cat{B}).
\]

The main insight of the proof is that a pair of trees as in
Figure~\ref{fig:twotrees} can be regarded as a span in the category of
forests, and multiplication of such pairs in the Thompson group is nothing
more than the usual composition of spans (by pullback).  The only significant
work in the proof is to establish the latter fact.

The universal property of $\cat{A}$ is immediate:
\begin{propn}
$\cat{A}$ is the monoidal category freely generated by an idempotent object.
\end{propn}

\begin{proof}
To lighten the notation, write $R_i (X)$ as $X$.  Then for any
monoidal category $\cat{M}$,
\begin{eqnarray}
\lefteqn{\MonCat (\cat{A}, \cat{M})}
\label{eq:beginning}	\\
	&\iso	&
\MonGpd (L_3 L_2 L_1 (\cat{B}), S_3 (\cat{M}))	
\\
	&\iso	&
\MonCat (L_2 L_1 (\cat{B}), S_3 (\cat{M}))	
\\
	&\iso	&
\Multicat (L_1 (\cat{B}), S_3 (\cat{M}))	
\\
	&\iso	&
\{ (M, G) 
\such 
M \in \cat{M}, \
G \in \Operad (L_1 (\cat{B}), \End_{S_3 (\cat{M})} (M)) \}	
\label{eq:mo}	\\
	&\iso	&
\{ (M, \mu) 
\such
M \in \cat{M}, \ 
\mu \in \Set^\nat (\cat{B}, \End_{S_3 (\cat{M})} (M) \}
\\
	&\iso	&
\{ \textrm{idempotent objects in } \cat{M} \}
\label{eq:Brep}
\end{eqnarray}
naturally in $\cat{M}$.  Most of these isomorphisms are by adjointness;
\bref{eq:mo}~is from the final observation in the section on multicategories
in the Appendix; \bref{eq:Brep}~is the fact that a map from $\cat{B}$ to
another signature $\cat{B'}$ just picks out an element of $\cat{B'}_2$, which
in this case is the set of maps $M^{\otimes 2} \go M$ in the groupoid
$S_3(\cat{M})$.

Hence $\cat{A}$ represents the functor $J: \MonCat \go \Set$ mapping a
monoidal category to the set of idempotent objects in it.  The generic
idempotent object $(A, \alpha) \in J(\cat{A})$ is obtained by tracing the
element $1_{\cat{A}}$ through the isomorphisms
\bref{eq:beginning}--\bref{eq:Brep}; then $(\cat{A}, A, \alpha)$ has the
universal property required.  
\done
\end{proof}

To obtain an explicit description of $(\cat{A}, A, \alpha)$, and
in particular of the automorphism group of $A$, we go through
each step of the construction using the descriptions of the adjoint functors
given in the Appendix.

First step: the free operad $\cat{C} = L_1(\cat{B})$ is the operad of
(unlabelled, binary) trees; thus, $\cat{C}_n = \bintr{n}$ and composition in
$\cat{C}$ is by gluing roots to leaves.

Second step: $\cat{D} = L_2 (\cat{C})$ is the monoidal category in which
objects are natural numbers and maps $n \go k$ are forests with $n$ leaves and
$k$ roots (as defined in~\S\ref{sec:term}).  Composition is by gluing; tensor
of objects is addition; tensor of maps is juxtaposition.

\begin{lemma}
The forest category $\cat{D}$ has pullbacks.  
\end{lemma}

\begin{proof}
Any map $T: n \go k$ in $\cat{D}$ decomposes uniquely as a tensor product
$T = T_1 \otimes \cdots \otimes T_k$ with $T_i: n_i \go 1$, so it suffices
to prove that every diagram of the form
\[
\begin{diagram}
m	&		&	&		&m'	\\
	&\rdTo<{(\tau)}	&	&\ldTo>{(\tau')}&	\\
	&		&1	&		&	\\
\end{diagram}
\]
has a pullback (where $\tau$ and $\tau'$ are trees with $m$ and $m'$ leaves
respectively).  Indeed, let $\rho$ be the smallest tree containing both $\tau$
and $\tau'$ as subtrees.  Then
\[
(\tau) \of (\sigma_1, \ldots, \sigma_m)
= 
\rho
=
(\tau') \of (\sigma'_1, \ldots, \sigma'_{m'})
\]
for unique $\sigma_i$ and $\sigma'_{i'}$.  Writing $p$ for the number of
leaves of $\rho$, the square
\[
\begin{diagram}
	&		&p	&		&	\\
	&\ldTo<{(\sigma_1, \ldots, \sigma_m)}
			&	&\rdTo>{(\sigma'_1, \ldots, \sigma'_{m'})}
						&	\\
m	&		&	&		&m'	\\
	&\rdTo<{(\tau)}	&	&\ldTo>{(\tau')}&	\\
	&		&1	&		&	\\
\end{diagram}
\]
is a pullback.  
\done
\end{proof}

Third step: $\cat{E} = L_3 (\cat{D})$ is the monoidal groupoid in which
objects are natural numbers and maps $k' \go k$ are equivalence classes of
spans
\[
\begin{diagram}
	&	&n	&	&	\\
	&\ldTo<{(\tau'_1, \ldots, \tau'_{k'})}
		&	&\rdTo>{(\tau_1, \ldots, \tau_k)}
				&	\\
k'	&	&	&	&k	\\
\end{diagram}
\]
in $\cat{D}$.  Equivalence is generated by declaring this span to be
equivalent to 
\[
\begin{diagram}
	&	&p	&	&	\\
	&
\ldTo<{(\tau'_1, \ldots, \tau'_{k'}) \of (\rho_1, \ldots, \rho_n)}
		&	&
\rdTo>{(\tau_1, \ldots, \tau_k) \of (\rho_1, \ldots, \rho_n)}
				&	\\
k'	&	&	&	&k	\\
\end{diagram}
\]
for any forest $(\rho_1, \ldots, \rho_n)$ with $n$ roots (writing $p$ for its
number of leaves), and it makes no difference if we insist that all but one of
the $\rho_i$s is trivial and the remaining one is the 2-leafed tree
$\binarytree$.

Final step: $\cat{A}$ is the underlying monoidal category of $\cat{E}$.  Under
the isomorphisms \bref{eq:beginning}--\bref{eq:Brep}, the identity $1_\cat{A}$
corresponds to the idempotent object $(1, \alpha)$ in $\cat{A}$, where
$\alpha$ is the equivalence class of the span
\[
\begin{diagram}
	&		&2	&                       &	\\
	&\ldTo<\id	&	&\rdTo>{(\binarytree)}	&	\\
2	&		&	&                       &1.	\\
\end{diagram}
\]
So to prove Theorem~\ref{thm:main}, we have to show that 
$\Aut_{\cat{A}}(1) \iso F$.

Since $\cat{A}$ is a groupoid, $\Aut_{\cat{A}}(1)$ consists of all maps $1
\go 1$ in $\cat{A}$.  We have just seen that such a map is an equivalence
class of pairs $(\tau, \tau')$ of trees with the same number of leaves,
where equivalence is generated by $[\tau, \tau'] = [\omega^n_i (\tau),
\omega^n_i (\tau')]$ whenever $\tau, \tau' \in \bintr{n}$ and $1 \leq i
\leq n$.  To compose maps
\[
1 \goby{[\sigma, \sigma']} 1 \goby{[\tau, \tau']} 1,
\]
form the diagram 
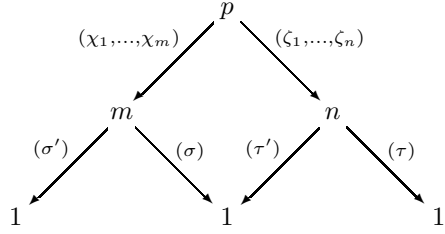
\[
\begin{diagram}
	&	&	&	&p	&	&	&	&	\\
	&	&	&\ldTo<{(\chi_1, \ldots, \chi_m)}
				&	&\rdTo>{(\zeta_1, \ldots, \zeta_n)}
						&	&	&	\\
	&	&m	&	&	&	&n	&	&	\\
	&\ldTo<{(\sigma')}
		&	&\rdTo>{(\sigma)}
				&	&\ldTo<{(\tau')}
						&	&\rdTo>{(\tau)}&\\
1	&	&	&	&1	&	&	&	&1	\\
\end{diagram}
\]
in which the square is a pullback; then
\[
[\tau, \tau'] \of [\sigma, \sigma']
=
[\tau \of (\zeta_1, \ldots, \zeta_n),\littlespace
\sigma' \of (\chi_1, \ldots, \chi_m)].
\]
There exist $i_1, \ldots, i_r, n_1, \ldots, n_r$ with the
property that for all $\pi \in \bintr{n}$,
\[
\pi \of (\zeta_1, \ldots, \zeta_n)
=
\omega_{i_r}^{n_r} \cdots \omega_{i_1}^{n_1} (\pi),
\]
and similarly $j_1, \ldots, j_s, m_1, \ldots, m_s$ for $(\chi_1,
\ldots, \chi_m)$.  Hence
\[
\omega^{n_r}_{i_r} \cdots \omega^{n_1}_{i_1} (\tau')
=
\omega^{m_s}_{j_s} \cdots \omega^{m_1}_{j_1} (\sigma)
\]
and 
\[
[\tau, \tau'] \of [\sigma, \sigma']
=
[\omega^{n_r}_{i_r} \cdots \omega^{n_1}_{i_1} (\tau),\littlespace
\omega^{m_s}_{j_s} \cdots \omega^{m_1}_{j_1} (\sigma')].
\]
But this description of $\Aut_{\cat{A}}(1)$ is exactly the description of
$F$ in~\S\ref{sec:term}.  Hence $\Aut_{\cat{A}}(1) \iso F$, proving
Theorem~\ref{thm:main}. 

Finally, we remark that the proof can be recast slightly so that the diagram
in Figure~\ref{fig:steps} becomes a chain of adjunctions
\[
\Set^\nat
\adjnpair
\Operad
\adjnpair
\Multicat_*
\adjnpair
\MonCat_*
\adjntriple
\MonGpd_*.
\]
Here $\Multicat_*$ denotes the category of multicategories equipped with a
distinguished object, and similarly $\MonCat_*$ and $\MonGpd_*$.  The content
of the argument is the same.

\appendix

\section{Appendix: Some categorical structures}

Here we review some categorical structures used in the proof of
Theorem~\ref{thm:main}: signatures, operads, multicategories, groupoids, and
monoidal groupoids.  We also review the basic relationships between these
structures.

\smallheading{Signatures}
We use the category $\Set^\nat$, the product of $\nat$ copies of the
category of sets.  Its objects are sequences $(\cat{B}_n)_{n\in\nat}$ of sets,
which can be regarded as \demph{signatures} for finitary, single-sorted
algebraic theories; $\cat{B}_n$ is thought of as the set of $n$-ary
operations.

\smallheading{Operads}
In this section and the next (on multicategories), trees will not be assumed
to be binary: any natural number of branches, including $0$, may grow out of
each vertex.

\begin{figure}
\[
%
\begin{array}{c}
\setlength{\unitlength}{1em}
\begin{picture}(16.5,12)(-0.5,-6)
\cell{10}{0}{l}{\tusual{\theta}}
\cell{2}{4}{l}{\tusual{\theta_1}}
\cell{2}{-4}{l}{\tusual{\theta_n}}
\cell{14}{0}{l}{\toutputrgt{C}}
\cell{2}{4}{r}{\tinputslft{C_1^1}{C_1^{k_1}}}
\cell{2}{-4}{r}{\tinputslft{C_n^1}{C_n^{k_n}}}
\qbezier(10,1.5)(8.5,1.5)(8,2.75)
\qbezier(6,4)(7.5,4)(8,2.75)
\cell{8.5}{2.75}{l}{C_1}
\qbezier(10,-1.5)(8.5,-1.5)(8,-2.75)
\qbezier(6,-4)(7.5,-4)(8,-2.75)
\cell{8.5}{-2.75}{l}{C_n}
\cell{9.2}{0.3}{c}{\vdots}
\cell{3}{0}{c}{\cdot}
\cell{3}{1}{c}{\cdot}
\cell{3}{-1}{c}{\cdot}
\end{picture}
\end{array}
\mbox{\hspace{1em}}
\textrm{gives}
\mbox{\hspace{2em}}
%
%
\begin{array}{c}
\setlength{\unitlength}{1em}
\begin{picture}(10,12)(0,-6)
\put(2,-6){\line(0,1){12}}
\put(8,0){\line(-1,-1){6}}
\put(8,0){\line(-1,1){6}}
\cell{4.9}{0}{c}{%
\scalebox{.9}{\ensuremath{\theta \of (\theta_1, \ldots, \theta_n)}}
}
\cell{2}{4}{r}{\tinputslft{C_1^1}{C_1^{k_1}}}
\cell{2}{-4}{r}{\tinputslft{C_n^1}{C_n^{k_n}}}
\cell{1.2}{0.3}{c}{\vdots}
\cell{8}{0}{l}{\toutputrgt{C}}
\end{picture}
\end{array}
\]
\caption{Composition in a multicategory}
\label{fig:opd-comp}
\end{figure}
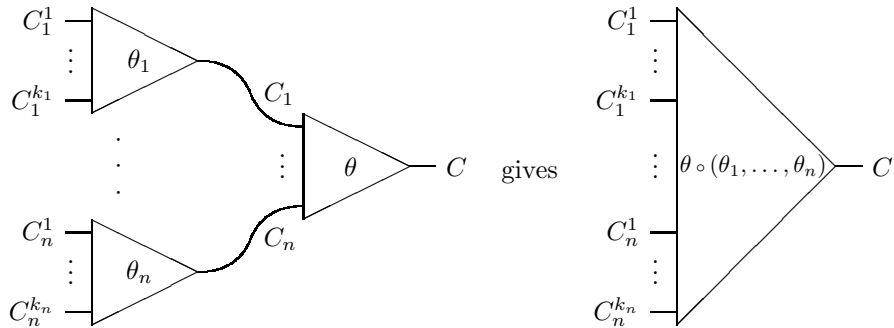

If $D$ is an object of a monoidal category $\cat{D}$ then the sequence
$(\cat{D}(D^{\otimes n}, D))_{n\in\nat}$ of hom-sets admits certain
algebraic operations.  This is the archetypal example of an operad.
Formally, an \demph{operad} consists of a sequence $(\cat{C}_n)_{n\in\nat}$
of sets together with a \demph{composition} function
\[
\begin{array}{ccc}
\cat{C}_n \times \cat{C}_{k_1} \times\cdots\times \cat{C}_{k_n}	&
\go		&
\cat{C}_{k_1 + \cdots + k_n}	\\
(\theta, \theta_1, \ldots, \theta_n)	&
\goesto		&
\theta \of (\theta_1, \ldots, \theta_n)
\end{array}
\]
for each $n, k_1, \ldots, k_n \in \nat$, and a \demph{unit} element $1 \in
\cat{C}_1$, satisfying associativity and unit axioms.  An element of
$\cat{C}_n$ can be thought of as an $n$-ary operation;
then composition is as shown in Figure~\ref{fig:opd-comp}
(ignoring the labels $C$, $C_i$, $C_i^j$).  The associativity and unit axioms
imply that every tree of operations has a
well-defined composite.

When $D$ is an object of a monoidal category $\cat{D}$, the `archetypal'
operad mentioned has composition
\begin{equation}	\label{eq:opd-comp}
\theta \of (\theta_1, \ldots, \theta_n)
=
\theta \of (\theta_1 \otimes\cdots\otimes \theta_n).
\end{equation}

(See~\cite{MSS}, \cite{May}, or~\cite{HOHC} for more on operads.  Monoidal
categories, operads and multicategories---see below---each come in symmetric
and non-symmetric versions.  We use the non-symmetric versions of everything.) 

There is an obvious notion of map of operads, giving a category $\Operad$.
Any operad has an underlying signature, giving a forgetful functor $\Operad
\go \Set^\nat$.  This has a left adjoint $L$; in other words, we may form the
free operad $L(\cat{B})$ on any signature $\cat{B}$.  The elements of
$(L(\cat{B}))_n$ are the $n$-leafed trees in which each vertex is
labelled by an element of $\cat{B}_k$, where $k$ is the number of branches
growing out of the vertex.  Composition in $L(\cat{B})$ is given
by gluing roots to leaves, and the unit is the trivial tree.  For the proof
and further details, see~\cite[2.3]{HOHC}.

\smallheading{Multicategories}
A \demph{multicategory} resembles a category in that it consists of objects,
arrows between objects, and a unique composite for every composable diagram
of arrows.  The only difference is the shape of the arrows, which in a
multicategory are of the form
\begin{equation}	\label{eq:multiarrow}
C_1, \ldots, C_n \goby{\theta} C 
\end{equation}
where $n \in \nat$ and $C_1, \ldots, C_n, C$ are objects.  Composition is
as shown in Figure~\ref{fig:opd-comp}; to each object $C$ there is assigned
an identity arrow $1_C: C \go C$; and associativity and identity axioms
hold, so that every tree of arrows has a well-defined composite.  The
details can be found in~\cite{Lam} or~\cite{HOHC}.

A typical example of a multicategory has vector spaces as objects and
multilinear maps as arrows.  Composition is given by~\bref{eq:opd-comp}.
Similarly, any monoidal category has an \demph{underlying multicategory}: the
objects are the same, and an arrow~\bref{eq:multiarrow} in the multicategory
is an arrow $C_1 \otimes\cdots\otimes C_n \go C$ in the monoidal category.

Let $\MonCat$ be the category of small monoidal categories and monoidal
functors (all strict, as usual).  Let $\Multicat$ be the category of small
multicategories and maps between them (defined in the obvious way).  We have
just defined a forgetful functor $\MonCat \go \Multicat$.  It has a left
adjoint $L$: given a multicategory $\cat{C}$, an object of $L(\cat{C})$ is a
finite sequence $(C_1, \ldots, C_n)$ of objects of $\cat{C}$, and an arrow in
$L (\cat{C})$ is a finite sequence of arrows in $\cat{C}$.  Thus, arrows
\[
C_1^1, \ldots, C_1^{k_1} \goby{\theta_1} C_1,
\quad
\ldots,
\quad
C_n^1, \ldots, C_n^{k_n} \goby{\theta_n} C_n
\]
in $\cat{C}$ give rise to an arrow 
\[
(C_1^1, \ldots, C_1^{k_1}, \ldots, C_n^1, \ldots, C_n^{k_n})
\goby{(\theta_1, \ldots, \theta_n)}
(C_1, \ldots, C_n)
\]
in $L (\cat{C})$.  Tensor product in $L (\cat{C})$ is concatenation of
sequences.  

For example, let $1$ be the terminal multicategory, which has one object and
one $n$-ary arrow for each $n \in \nat$.  Then $L(1)$ is $\fcat{FinOrd}$, the
monoidal category of finite ordinals and order-preserving maps.

A multicategory $\cat{C}$ with only one object is just as an operad: if we
call the object $C$ and write $\cat{C}_n$ for the set of arrows
\begin{equation}	\label{eq:nfoldendo}
\underbrace{C, \ldots, C}_n \go C
\end{equation}
then the multicategory structure on $\cat{C}$ is exactly an operad
structure on $(\cat{C}_n)_{n\in\nat}$.  We write this operad as $\cat{C}$,
too. 

More generally, every object $C$ of a multicategory $\cat{C}$ has an
\demph{endomorphism operad} $\End_{\cat{C}}(C)$, whose $n$-ary operations are
the maps~\bref{eq:nfoldendo} in $\cat{C}$.

There is a full and faithful inclusion functor $\Operad \rIncl \Multicat$.
If $\cat{C}'$ is an operad and $\cat{C}$ a multicategory then a
map $\cat{C}' \go \cat{C}$ of multicategories amounts to an object $C \in
\cat{C}$ together with a map $\cat{C}' \go \End_{\cat{C}}(C)$ of operads.

\smallheading{Groupoids}
A \demph{groupoid} is a category in which every map is an
isomorphism.  The inclusion $\Gpd \rIncl \Cat$ from small
groupoids into small categories has a left adjoint, `free
groupoid'.  
It is slightly tricky to describe the free groupoid on an arbitrary
category~\cite{Pare2}, but it is straightforward when the category has
pullbacks~\cite{Pare1}.  This case is all that we will need.

We use the notion of bicategory~\cite{Ben}.  Given a category $\cat{D}$ with
pullbacks, first form the bicategory $\Span(\cat{D})$ of spans in
$\cat{D}$~\cite[2.6]{Ben}.  Then form the groupoid $\cat{E}$ whose objects are
those of $\cat{D}$, whose maps $D' \go D$ are the connected-components of the
hom-category $(\Span(\cat{D}))(D', D)$, and whose composition is inherited
from $\cat{D}$.  (This is made possible by the fact that the
connected-components functor $\Cat \go \Set$ preserves finite products.)

The groupoid $\cat{E}$ can be described explicitly.  A \demph{span} from $D'$
to $D$ is a diagram in $\cat{D}$ of the form
\[
\begin{diagram}
	&		&X	&		&	\\
	&\ldTo<{\phi'}	&	&\rdTo>\phi	&	\\
D'	&		&	&		&D,	\\
\end{diagram}
\]
written $(\phi, \phi')$; note the reversal.  Call two such spans $(\phi,
\phi')$, $(\psi, \psi')$ \demph{equivalent} if there exists a commutative
diagram of the form
\[
\begin{diagram}
	&		&X	&		&	\\
	&\ldTo<{\phi'}	&\uTo	&\rdTo>\phi	&	\\
D'	&		&Z	&		&D,	\\
	&\luTo<{\psi'}	&\dTo	&\ruTo>\psi	&	\\
	&		&Y	&		&	\\
\end{diagram}
\]
and write $[\phi, \phi']$ for the equivalence class of a span $(\phi,
\phi')$.  Then the objects of the groupoid $\cat{E}$ are those of $\cat{D}$,
the maps from $D'$ to $D$ in $\cat{E}$ are the equivalence classes of spans
from $D'$ to $D$ in $\cat{D}$, composition is by pullback, and $[\phi,
\phi']^{-1} = [\phi', \phi]$.

There is a functor $\eta_{\cat{D}}: \cat{D} \go \cat{E}$ defined on objects
as the identity and on maps by
\[
\left(
D' \goby{\theta} D
\right)
\quad
\goesto 
\quad 
[\theta, 1_{D'}].
\]
It is straightforward to check that every functor from $\cat{D}$ to a groupoid
factors uniquely through $\eta_{\cat{D}}$; hence $\cat{E}$ is the free
groupoid on $\cat{D}$.

Let $\Cat_\pull$ be the full subcategory of $\Cat$ consisting of
the categories with pullbacks.  Since every groupoid has
pullbacks, there is a forgetful functor $R: \Gpd \go \Cat_\pull$,
and we have just constructed its left adjoint $L$.  It is clear from the
construction that $L$, as well as $R$, preserves finite products.

\smallheading{Monoidal groupoids}
Let $\MonGpd$ be the full subcategory of $\MonCat$ consisting of the
monoidal groupoids.  The inclusion $\MonGpd \rIncl \MonCat$ also has a left
adjoint, which again is easily described in the presence of pullbacks.

Let $\MonCat_\pull$ be the full subcategory of $\MonCat$ consisting of the
monoidal categories with pullbacks.  Taking internal monoids throughout the
adjunction $\Gpd \oppair{R}{L} \Cat_\pull$ gives an adjunction $\MonGpd
\oppair{R}{L} \MonCat_\pull$.  This new $R$ is the evident forgetful functor.
If $\cat{D} \in \MonCat_\pull$ then $L(\cat{D})$ is the free groupoid on
$\cat{D}$ as constructed above, with monoidal structure inherited from
$\cat{D}$ in the obvious way.

The forgetful functor $\MonGpd \go \MonCat$ also has a
\emph{right} adjoint, sending a monoidal category to its
underlying groupoid (the subcategory consisting of all the
objects and all the isomorphisms).

\end{document}